\documentclass[11pt]{amsart}
\usepackage{amsfonts,amssymb,epsfig}
\usepackage[latin1]{inputenc}
\usepackage[all]{xy}
\usepackage{amsmath}
\usepackage{delarray}
\usepackage{amssymb}
\newtheorem{thm}{Theorem}[section]

\newtheorem{lem}[thm]{Lemma}

\theoremstyle{definition}

\theoremstyle{remark}

\newtheorem{remark}[thm]{Remark}
\numberwithin{equation}{section}

\renewcommand{\epsilon}{\varepsilon}

\newcommand{\cO}{\mathcal O}

\newcommand{\cU}{\mathcal U}
\newcommand{\bbR}{\mathbb R}

\newcommand{\bbN}{\mathbb N}

\begin{document}

\title
{Linearization of proper groupoids}

\author{Nguyen Tien Zung}
\address{Labo. Emile Picard, UMR 5580 CNRS, UFR MIG, Université Toulouse III}
\email{tienzung@picard.ups-tlse.fr} 

\date{Second version, May/2003}

\begin{abstract} We prove the following result, conjectured by Alan Weinstein:
every smooth proper Lie groupoid near a fixed point is locally
linearizable, i.e. it is locally isomorphic to the associated groupoid of
a linear action of a compact Lie group. In combination with a slice
theorem of Weinstein, our result implies the smooth linearizability of a
proper Lie groupoid in the neighborhood of an orbit under a mild
condition.
\end{abstract}
\maketitle

\section{Introduction}

In the theory of Lie groupoids, the role played by proper Lie groupoids is probably
comparable to that played by compact Lie groups in the Lie group theory. However,
while compact Lie groups have been classified completely long time ago, the
study of proper Lie groupoids has been initiated only recently by Alan
Weinstein \cite{Weinstein-Linearization2000,Weinstein-RPG2001}, and many basic
questions remained open. One of them, formulated by Weinstein as the \emph{proper
groupoid structure conjecture}, is about the local linearizability of proper Lie
groupoids near fixed points. The main purpose of the present paper is to prove this
conjecture.

So let us consider a smooth proper groupoid near a fixed point. It means that we are
given the following data:

i) A smooth manifold $M$ (the \emph{total space}), a smooth submanifold $B$ of $M$
(the \emph{base space}, also called the \emph{space of units}), and two smooth maps
$M \stackrel{s}{\longrightarrow} B$, $M \stackrel{t}{\longrightarrow} B$, called the
\emph{source map} and the \emph{target map} respectively. The restrictions of $s$
and $t$ on $B$ are identity. Here, by \emph{smooth} we mean either
$C^{\infty}$-smooth, or $C^{h}$-smooth for some natural number $h$.

ii) Denote by $M_k$ $(k \geq 2)$ the subset of $M \times \dots \times M$ ($k$ times)
consisting of $k$-tuples $(p_1,...,p_k)$  of elements of $M$ such that $s(p_1) =
t(p_2), s(p_2) = t(p_3), ..., s(p_{k-1}) = t(p_k)$. Then there is a map from $M_2$
to $M$, called the \emph{product} map: if $(p,q) \in M_2$, i.e. $s(p) = t(q)$, then
we can form the product of $p$ with $q$, denoted by $p.q$ .

iii) There is a map from $M$ to itself, called the \emph{inversion}, which
associates to each element $p \in M$ an element denoted by $p^{-1}$,
called its inverse. The inverse of  each point in $B$ is itself.

The product map and the inversion map satisfy the following conditions:
$s(p.q) = s(q)$ and $t(p.q) = t(p)$ for any $(p,q) \in M_2$; $(p.q).r =
p.(q.r)$ for any $(p,q,r) \in M_3$; $p.p^{-1} \in B$ and $p^{-1}.p \in B$
for any $p \in M$; $p.q = p$ if $q \in B$ and $p.q = q$ if $q \in B$.

iv) A point $x_0 \in B$, which is a fixed point for the groupoid, in the sense that
for any $p \in M$ with $s(p) = x_0$ we also have $t(p) = x_0$. The set $G = \{p \in
M, s(p) = x_0 \}$ is then a group under the above product map.

The properness conditions that we impose on $G$ and $s$ and $t$ are that $G$ is a
compact (not necessarily connected) Lie group, and $M \stackrel{s}{\longrightarrow}
B$, $M \stackrel{t}{\longrightarrow} B$ are smooth fiber bundles, with fibers
diffeomorphic to $G = s^{-1}(x_0) = t^{-1}(x_0)$.

For each point $x \in B$, the \emph{orbit} of our groupoid via $x$ is the
set $O(x) := t(s^{-1}(x))$, which is the same as $s(t^{-1}(x))$. A subset
in $B$ is called \emph{invariant} if it is saturated by the orbits of the
groupoid. As observed by Weinstein \cite{Weinstein-RPG2001}, any
neighborhood of $x_0$ in $B$ contains an invariant neighborhood (simply
take a sufficiently small open ball $D \subset B$ which contains $x_0$ and
take $D' = t(s^{-1}(D))$, then $D' \supset D$ is a small invariant open
neighborhood of $x_0$). Of course, the restriction of a groupoid to an
invariant set is still a groupoid. Since we are interested in the local
behavior of our proper groupoid $M \rightrightarrows B$ near $G =
s^{-1}(x_0)$, we can assume that $B$ is a small neighborhood of $x_0$ in
an Euclidean space.

If the compact Lie group $G$ acts smoothly on a manifold $X$, then we can
define the associated \emph{action groupoid} $G \times X \rightrightarrows
X$: the source map is $s(g, x) = x $ for $(g,x) \in G \times X$, the
target map is $t(g,x) = g.x$ (the action of $g$ on $x$), and the inversion
map is $(g,x)^{-1} = (g^{-1},g.x)$. The space of units is $\{1_G\} \times
X$ (where $1_G$ denotes the neutral element of $G$), which is identified
with $X$. If the action of $G$ on $X$ has a fixed point, then due to
Bochner's linearization theorem \cite{Bochner-Linearization1945}, we may
assume that $X$ is a small neighborhood of this fixed point in an
Euclidean space, and the action of $G$ on $X$ is linear. In this case we
say that the action groupoid $G \times X \rightrightarrows X$ is a
\emph{linear action groupoid}. We will say that our proper groupoid $M
\rightrightarrows B$ is isomorphic to an action groupoid $G \times X
\rightrightarrows X$ if there is a smooth diffeomorphism from $M$ to $G
\times X$ which preserves the maps $s,t$, the product and the inversion.

The main result of the present paper is the following

\begin{thm}
\label{thm:LPG} Let $M \rightrightarrows B$ be a $C^h$ (resp., $C^\infty$) smooth
proper groupoid with a fixed point $x_0 \in B$, with the notations as of above. Then
the restriction of this groupoid to a sufficiently small invariant neighborhood of
$x_0$ in $B$ is $C^h$ (resp., $C^\infty$) isomorphic to a linear action groupoid of
$G$ .
\end{thm}

An important application of Theorem \ref{thm:LPG} is in the study of
symplectic groupoids, where one can show that proper symplectic groupoids
are locally isomorphic to standard models, see
\cite{Zung-PGPC2003,Zung-TorusSurvey2003}. In turn, this result about the
local structure of proper symplectic groupoids is important for the study
of intrinsic convexity properties of momentum maps, initiated by Weinstein
\cite{Weinstein-Linearization2000,Weinstein-Convexity2001}, see
\cite{Zung-PGPC2003,Zung-TorusSurvey2003}.

Theorem \ref{thm:LPG}, together with Weinstein's slice theorem (Theorem
9.1 of \cite{Weinstein-RPG2001}), immediately implies the following result
about the local linearization of a proper groupoid in the neighborhood of
an orbit:

\begin{thm}
\label{thm:LPG-Orbit} Let $\Gamma \rightrightarrows X$ be a proper Lie
groupoid, and let $\cO$ be an orbit of $\Gamma$ which is a manifold of
finite type. Then there is a neighborhood $\cU$ of $\cO$ in $X$ such that
the restriction of $\Gamma$ to $\cU$ is isomorphic to the restriction of
the action groupoid $\Gamma_\cO \times_\cO N\cO \rightrightarrows N\cO$ to
a neighborhood of the zero section.
\end{thm}

In Theorem \ref{thm:LPG-Orbit}, $N\cO$ means the normal bundle of $\cO$ in
$X$, and $\Gamma_\cO$ is the restriction of $\Gamma$ to $\cO$. The
condition that $\cO$ is of finite type means that there is a proper map
from $\cO$ to $\bbR$ with a finite number of critical points. See
\cite{Weinstein-RPG2001} for details. \\

Theorem \ref{thm:LPG} is essentially equivalent to the existence of a
smooth \emph{surjective homomorphism $\phi$ from $M$ to $G$} (after
shrinking $B$ to a sufficiently small invariant neighborhood of $x_0$),
i.e. a smooth map $\phi: M \to G$ which satisfies
\begin{equation}
\label{eqn:homomorphism} \phi(p.q) = \phi(p). \phi(q) \ \ \forall \ (p,q) \in M_2 \
,
\end{equation}
such that the restriction of $\phi$ to $G = s^{-1}(x_0) \subset M$ is an
automorphism of $G$. We may assume that this automorphism is identity.

Indeed, if there is an isomorphism from $M$ to an action groupoid with the total
space $G \times X$, then the composition of the isomorphism map $M \to G \times X$
with the projection $G \times X \to G$ is such a homomorphism. Conversely, if we
have a homomorphism $\phi: M \to G$, whose restriction to $G = s^{-1}(x_0) \in M$ is
the identity map of $G$, then assuming that $B$ is a sufficiently small invariant
neighborhood of $x_0$, we have a diffeomorphism
\begin{equation}
(\phi, s): M \to G \times B \ .
\end{equation}
Denote by $\theta$ the inverse map of $(\phi, s)$.
Then there is an action of $G$ on $B$ defined by $g.x = t(\theta(g,x))$, and the map
$(\phi, s)$ will be an isomorphism from $M \rightrightarrows B$ to the action
groupoid $G \times B \rightrightarrows B$. This action groupoid is linearizable by
Bochner's theorem, implying that the groupoid $M \rightrightarrows B$ is
linearizable.

In order to find such a homomorphism from $M$ to $G$, we will use the averaging
method. The idea is to start from an arbitrary smooth map $\phi: M \to G$ such that
$\phi|_G = Id$ (Recall that we identify $s^{-1}(x_0) = t^{-1}(x_0)$ with $G$). Then
Equality (\ref{eqn:homomorphism}) is not satisfied in general, but it is satisfied
for $p, q \in G$. Hence it is ``nearly satisfied'' in a small neighborhood of $G =
s^{-1}(x_0)$ in $M$. In other words, if the base $B$ is small enough, then
$\phi(p.q)\phi(q)^{-1}$ is near $\phi(p)$ for any $(p,q) \in M_2$. We will replace
$\phi(p)$ by the average value of $\phi(p.q)\phi(q)^{-1}$ for $q$ running on
$t^{-1}(s(p))$ (it is to be made precise how to define this average value). This way
we obtain a new map $\widehat{\phi}: M \to G$, which will be shown to be ``closer''
to a homomorphism than the original map $\phi$. By iterating the process and taking
the limit, we will obtain a true homomorphism $\phi_{\infty}$ from $M$ to $G$.

Let us mention that the averaging method is a classical method which is relatively
simple and effective in many problems involving compact group actions. For example,
it was used by Elie Cartan \cite{Cartan-Riemann1928,Cartan-Riemann1929} to prove the
existence of a fixed point of a compact group action under some conditions
(see \cite{Weinstein-Invariant2000} and references therein for recent results on
invariant manifolds of compact group actions).
Bochner's linearization theorem
\cite{Bochner-Linearization1945} is indeed based on a very simple averaging formula.
The iterative averaging method for was already used by Grove, Karcher and Ruh in
\cite{GrKaRu-Group1974} to show that near-homomorphisms between compact Lie groups
can be approximated by homomorphisms. Weinstein
\cite{Weinstein-Invariant2000,Weinstein-Linearization2000,Weinstein-RPG2001} also
suggested the averaging method for proving Theorem \ref{thm:LPG}, although he
considered near-homomorphisms from $G$ to the group of bisections of $M
\rightrightarrows B$ instead of near-homomorphisms from $M$ to $G$.
In this paper, instead of using the Haar
measure of a compact group for making the averaging, we will use a ``Haar system''
on a proper groupoid. More precisely, we will use a family of probability
measures on the fibers of the fibration $M \stackrel{t}{\longrightarrow} B$, which
are invariant under left translations. This ``technical'' detail allows us to
effectively control the convergence of our iterative averaging process, using
standard Banach-space estimations.

It is clear (see, e.g., \cite{Weinstein-Linearization2000}) that the problem of
linearization of Lie groupoids is closely related to the problem of linearization of
Lie algebroids. In this aspect, Theorem \ref{thm:LPG} is related to our results with
Philippe Monnier \cite{Zung-Levi2002,MoZu-Levi2002} on Levi decomposition of Lie
algebroids. Several natural questions arise, among which: is there a similar result
about Levi decomposition of Lie groupoids near a fixed point ? If it exists, can we
prove it by the averaging method ? Let us mention that the results of
\cite{Zung-Levi2002,MoZu-Levi2002} are obtained with the aid of the fast convergence
methods of Kolmogorov and Nash-Moser, which are technically considerably more
complicated than the averaging method. Theorem \ref{thm:LPG} or a more general
result about Levi decomposition of Lie groupoids, together with criteria about
(partial) integrability of Lie algebroids like the one of Crainic-Fernandes
\cite{CrFe-Lie2001}, might eventually lead to a simplification and improvement of
our papers \cite{Zung-Levi2002,MoZu-Levi2002}.

The rest of this paper is devoted to the construction of a smooth homomorphism from
$M$ to $G$ (after shrinking the base $B$ if necessary), whose restriction to $G$ is
the identity map of $G$. As explained above, Theorem \ref{thm:LPG} follows
immediately from the existence of such a homomorphism. In Section 2 we present our
iterative averaging algorithm, and in Section 3 we use standard analytical estimates
to show that our algorithm actually yields a ($C^h$ or $C^\infty$) smooth
homomorphism.

In this paper, we will consider only $C^h$ and $C^\infty$ smooth groupoids,
but we suspect that, with more care, one can probably show that our algorithm for
constructing a linearization of $M \rightrightarrows B$ works in the real-analytic
category as well.

\section{The averaging algorithm}

\subsection{A translation-invariant measure}

For each point $p \in M$ denote by $T(p) = t^{-1}(t(p))$ (resp., $S(p) =
s^{-1}(s(p))$) the fiber of the fibration $M \stackrel{t}{\longrightarrow} B$ (resp.
$M \stackrel{s}{\longrightarrow} B$) which contains $p$.

Each element $q \in M$ defines a {\it left translation map} $\tau_q$ from $T(s(q))$
to $T(q)$ as follows: $\tau_q(r) = q.r$ for $r \in T(s(q))$.

Let us fix on $M$ an arbitrary smooth Riemannian metric. (Throughout this paper,
smooth objects have the same smoothness class as that of our groupoid $M
\rightrightarrows P$, unless mentioned explicitly otherwise). It induces on each
fiber $T(p)$ of the fibration $M \stackrel{t}{\longrightarrow} B$ a smooth measure
(which depends smoothly on the fiber itself), which we will denote by $d\mu_0$.
Similarly, the metric on $M$ induces on each fiber $S(p)$ of the fibration $M
\stackrel{s}{\longrightarrow} B$ a smooth measure, which we will denote by $d\nu_0$.

\begin{lem}
There is a smooth positive function $f$ on $M$ such that the measure $d\mu =
fd\mu_0$ defined on each fiber of the fibration $M \stackrel{t}{\longrightarrow} B$
is a translation-invariant probability measure: for each $q \in M$ the left
translation map $\tau_q$ from $T(s(q))$ to $T(q)$ preserves the measure $d\mu$, and
moreover $\int_{T(q)}1 d\mu = 1$.
\end{lem}

Translation-invariant measures are also called {\it Haar
systems} (see, e.g., \cite{AnRe-Groupoid2000}), and the above lemma is
probably a well-known folklore result , but to make our paper
self-contained, we will include a simple proof of it here.

{\it Proof}. Define $\widetilde{f}$ by the following integral formula:

\begin{equation}
\widetilde{f}(r) = \int_{p \in S(r)}
\frac{\tau_{p.r^{-1}}^{\ast}d\mu_0(r)}{d\mu_0(r)} d\nu_0 \ .
\end{equation}

In the above formula, $\tau_{p.r^{-1}}$ is the left translation map from $T(r)$ to
$T(p)$, $\tau_{p.r^{-1}}^{\ast}d\mu_0(r)$ means the pull-back at $r$ of the smooth
measure $d\mu_0$ on $T(p)$ by $\tau_{p.r^{-1}}$, and $d\nu_0$ is the measure on
$S(r)$ induced by the metric on $M$. One checks immediately that $\widetilde{f}$ is
a smooth positive function on $M$, and the measure $d\widetilde{\mu} =
\widetilde{f}d\mu_0$ on the fibers of $M \stackrel{t}{\longrightarrow} B$ is
translation invariant. Dividing $\widetilde{f}$ by the integral $I = \int_{T(p)}
\widetilde{f}d\mu_0$, we get the required smooth function $f = \widetilde{f}/I$ and
translation-invariant measure $d\mu = f d\mu_0$. \hfill $\square$

From now on, we will fix a measure $d\mu$ provided by the above lemma. We will use
it in our averaging formulas. The invariance of $d\mu$ under translations will be
important for us.

\begin{remark} As pointed out to me by A. Weinstein, there is a
natural 1-1 correspondence between smooth translation-invariant volume
forms  on fibers of $M \stackrel{t}{\rightarrow} B$ and smooth
non-vanishing sections of the top exterior power of the associated Lie
algebroid over $B$. This point of view gives another simple proof of the
above lemma, without the need of choosing a Riemannian metric on $M$.
\end{remark}

\begin{remark} Another simple way to define a measure on $T(q)$ is to
restrict the map $\phi: M \to G$ to $T(q)$ and take the pull-back of the
Haar measure of $G$ via this restricted map. But this measure is not
translation-invariant, and moreover it depends on the differential of
$\phi$. Maybe the averaging algorithm presented in the next subsection
still works with this induced measure from $G$, but the estimates involved
will be much more complicated, and one might have to use Nash-Moser theory
due to ``loss of differentiability'' in some inequalities. (This loss is
hidden in the measure, which depends on the differential of $\phi$).
\end{remark}

\subsection{The averaging formula}

We fix a bi-invariant metric on the Lie algebra $\mathfrak{g}$ of $G$ and the
induced bi-invariant metric on $G$ itself. Denote by $1_G$ the neutral element of
$G$. For each number $\rho > 0$, denote by $B_{\mathfrak{g}}(\rho)$ (resp.,
$B_{G}(\rho)$) the closed ball of radius $\rho$ in $\mathfrak{g}$ (resp., $G$)
centered at $0$ (resp., $1_G$). By rescaling the metric if necessary, we will assume
that the exponential map
\begin{equation}
\exp : B_{\mathfrak{g}}(1) \to B_{G}(1)
\end{equation}
is a diffeomorphism. Denote by
\begin{equation}
\log : B_{G}(1) \to B_{\mathfrak{g}}(1)
\end{equation}
the inverse of $\exp$. Define the distance $\Delta(\phi)$ of $\phi$ from being a
homomorphism as follows:

\begin{equation}
\Delta(\phi) = \sup_{(p,q) \in M_2} d(\phi(p.q).\phi(q)^{-1}.\phi(p)^{-1}, 1_G) \ .
\end{equation}

Let $\phi: M \to G$ be a smooth map such that $\phi|_G$ is identity.  We will assume
that $\Delta(\phi) \leq 1$, so that the following map $\widehat{\phi}: M \to G$ is
clearly well-defined:

\begin{equation}
\widehat{\phi}(p) = \exp( \int_{q \in T(s(p))}
\log(\phi(p.q).\phi(q)^{-1}.\phi(p)^{-1}) d\mu). \phi(p) \ .
\end{equation}

Since $d\mu$ is translation-invariant, by the change of variable $r = p.q$, we can
also write $\widehat{\phi}$ as:

\begin{equation}
\widehat{\phi}(p) = \exp( \int_{r \in T(p)}
\log(\phi(r).\phi(p^{-1}.r)^{-1}.\phi(p)^{-1}) d\mu). \phi(p) \ .
\end{equation}

Due to the commutativity of the maps $\exp$ and $\log$ with the adjoint actions, we
can also write $\widehat{\phi}$ as:

\begin{equation}
\widehat{\phi}(p) = \phi(p). \exp( \int_{q \in T(s(p))}
\log(\phi(p)^{-1}.\phi(p.q).\phi(q)^{-1}) d\mu) \ .
\end{equation}

It is clear that $\widehat{\phi}$ is a smooth map from $M$ to $G$, and its
restriction to $G = s^{-1}(0) \subset M$ is also identity.

\subsection{The iterative process}

Starting from an arbitrary $\phi$ (such that $\phi|_G = Id$) as above, we construct
a sequence of maps $\phi_n: M \to G$, by the recurrence formula $\phi_1= \phi$,
$\phi_{n+1} = \widehat{\phi_n}$. In the next section we will show that this sequence
is well-defined (after a restriction of $B$ to a smaller invariant neighborhood of
$x_0$ if necessary), and that
\begin{equation}
\phi_\infty = \lim_{n \to \infty} \phi_n
\end{equation}
exists, is smooth, and is a homomorphism from $M$ to $G$.

\section{Proof of convergence}

\subsection{Spaces of maps and $C^k$-norms}

As in the introduction, for each $n \in \bbN,$ we denote by $M_n$ the space of
$n$-tuples of points $(p_1,...,p_n)$ of $M$ such that $s(p_1)= t(p_2),
s(p_2)=t(p_3),\dots, s(p_{n-1}) = t(p_n)$. There are many natural product and
projection maps among these spaces $M_n$ ($n \in \bbN$), and these maps are smooth.
The reader may recall that these spaces $M_n$, together with natural maps among
them, are used in the definition of the cohomology of $M$, though we will not use
any cohomology (at least not in an explicit way) here. For each $n \in \bbN$, the
manifold $M_n$ is smoothly diffeomorphic to $B \times G \times ... \times G$ ($n$
copies of $G$). To fix the norms, we will fix such a diffeomorphism for each $n$. We
will mainly use the manifolds $M (= M_1), M_2$ and $M_3$. To fix the norms on $B$
(i.e. for maps from and to $B$), we will assume that $B$ is a neighborhood of $x_0$
in a given Euclidean space. (We will shrink $B$ whenever necessary, but the norm of
the Euclidean space which contains it will not be changed).

If $V_1$ and $V_2$ are two nonnegative numbers which depend on several variables and
parameters, then we will write $V_1 \preceq V_2$ (read $V_1$ is smaller than $V_2$
up to a multiplicative constant) if there is a positive constant $C$ (which does not
depend on the movable parameters and variables of $V_1$ and $V_2$) such that $V_1
\leq C V_2$. We can also write $V_1 = O(V_2)$ using Landau notation. We will write
$V_1 \approx V_2$ if $V_1 \preceq V_2 \preceq V_1$.

We are interested in the $C^k$-topology ($k \in {\mathbb Z}_+$, $k \leq h$ if $M$ is
only $C^h$-smooth) of the spaces of maps from $M,M_2,M_3$ to $\mathfrak{g}$ and $G$.
We will use $\|.\|_k$ to denote a fixed $C^k$-norm. (It doesn't matter much which
$C^k$-norm we use because they are equivalent). We use the following $C^0$ norms for
functions from $N$ (where $N$ denotes one of the spaces $M, M_2,M_3,...$) to $G$:
\begin{equation}
\|f\|_0 = \sup_{x \in N} d(f(x), 1_G) \ ,
\end{equation}
where $d(.,.)$ is the metric on $G$. To define a $C^k$-norm for functions from $N$
to $G$, we identify $TG$ with $\mathfrak{g} \times G$ and so on. Our convention is
that the constant map from $N$ to the neutral element in $G$ is also denoted by
$1_G$, and $\|1_G\|_k = 0 \ \forall k$. Moreover, for each fixed $k$, if $f$ is a
map from $N$ to $\mathfrak{g}$ with $\|f\|_0 \preceq 1$ and $\exp(f)$ is the
corresponding map from $N$ to $G$ then

\begin{equation}
\|f\|_k \approx \|\exp(f)\|_k \ .
\end{equation}

Let us write down some standard inequalities which will be used later on.

If $m$ is a fixed smooth map from $N$ to $N'$ (e.g., the product map from $M_2$ to
$M$) and $f$ is a map from $N'$ to $\mathfrak{g}$ or $G$, then we have (for each
fixed nonnegative integer $k$ which does not exceed the smoothness class of the
groupoid):

\begin{equation}
\label{eqn:composition_est} \|f\circ m \|_k \preceq \|f\|_k \ .
\end{equation}

If $f_r$ from $N$ to $\mathfrak{g}$ depends on a parameter $r$ which lives in a
probability space $R$, then

\begin{equation}
\label{eqn:Ck_integral} \|\int_R f_r dr \|_k \leq \sup_{r \in R} \|f_r\|_k \ .
\end{equation}

If $f_1, f_2$ are two functions from $N$ to $G$ then we have:
\begin{equation}
\label{eqn:C0_product} \| f_1.f_2 \|_0 \preceq  \|f_1\|_0 + \|f_2\|_0  \ ,
\end{equation}
and if moreover $\|f_1\|_{k-1}, \|f_2\|_{k-1} \preceq 1$ (for some fixed $k \geq
1$), then we have:
\begin{equation}
\label{eqn:Ck_product} \| f_1.f_2 \|_k \preceq  \|f_1\|_k + \|f_2\|_k
\end{equation}
and (more refined inequalities)
\begin{equation}
\label{eqn:Ck_product2} \| f_1.f_2 \|_k - \|f_2\|_k \preceq  \|f_1\|_k +
\|f_1\|_0\|f_2\|_k  \ ,
\end{equation}
\begin{equation}
\label{eqn:Ck_conjugate} \| f_1.f_2.f_1^{-1} \|_k \preceq  \|f_2\|_k + \|f_1\|_k.
\|f_2\|_0  \ .
\end{equation}

 For $f_1, f_2: N \to \mathfrak{g}$, define  $\log (\exp(f_1).\exp(f_2))$ by the
Baker-Campbell-Haussdorff formula. Then we have

\begin{equation}
\label{eqn:C0-exp} \|\log(\exp(f_1).\exp(f_2)) - f_1 - f_2 \|_0 \preceq
\|f_1\|_0\|f_2\|_0 \ ,
\end{equation}
and
\begin{equation}
\label{eqn:C0-exp2} \|\exp(f_1).\exp(f_2)\|_0 \preceq \|f_1+f_2\|_0 \ .
\end{equation}

If, moreover, $\|f_1\|_{k-1}, \|f_2\|_{k-1} \preceq 1$ (for some fixed $k \geq 1$)
then

\begin{equation}
\label{eqn:Ck-exp} \|\log(\exp(f_1).\exp(f_2)) - f_1 - f_2 \|_k \preceq
\|f_1\|_0\|f_2\|_k + \|f_2\|_0\|f_1\|_k + \|f_1\|_{k-1} +  \|f_2\|_{k-1} \ ,
\end{equation}
and
\begin{equation}
\label{eqn:Ck-exp2} \|\exp(f_1).\exp(f_2)\|_k \preceq \|f_1+f_2\|_k +
\|f_1\|_0\|f_2\|_k + \|f_2\|_0\|f_1\|_k + \|f_1\|_{k-1} +  \|f_2\|_{k-1} \ .
\end{equation}

Finally, if $f_n$ ($n\in \bbN$) are maps from $N$ to $G$, and $a_n$ are positive
numbers such that $\sum_{n=1}^\infty a_n$ converges and $\|f_n\|_k \leq a_k$ (for
some nonnegative integer $k$), then the product $f_n.f_{n-1}...f_{1}$ converges in
$C^k$-topology when $n \to \infty$ to a $C^k$-map from $N$ to $G$.

\subsection{$C^0$ estimates}

\begin{lem}
\label{lem:C0-estimate} For any $\phi: M \to G$ with $\Delta(\phi) \leq 1$ we have
\begin{equation}
\Delta(\widehat{\phi}) \preceq (\Delta(\phi))^2 \ .
\end{equation}
In particular, there is a positive constant $C_0 > 0, C_0 \leq 1$ such that if $\Delta(\phi)
\leq C_0$ then $\widehat{\phi}$ is well-defined and
\begin{equation}
\Delta(\widehat{\phi}) \leq (\Delta(\phi))^2 / C_0 \leq \Delta(\phi) \ .
\end{equation}
\end{lem}

Proof: Denote

\begin{equation}
\psi(p,q) = \phi(p.q).\phi(q)^{-1}. \phi(p)^{-1} \ ,
\end{equation}
and
\begin{equation}
\widehat{\psi}(p,q) =
\widehat{\phi}(p.q).\widehat{\phi}(q)^{-1}.\widehat{\phi}(p)^{-1} \ .
\end{equation}

Then $\psi$ and $\widehat{\psi}$ are functions from $M_2$ to $G$. By definition of
$\widehat{\phi}$, we have

\begin{equation}
\begin{array}{lll}
 \widehat{\psi}(p,q)  & =
  & \exp(\int_{r \in T(s(q))} \log(\psi(p.q,r)) d\mu). \phi(p.q). \phi(q)^{-1}. \\
& & \exp(\int_{r \in T(s(q))}^{-1} \log(\psi(q,r)) d\mu)^{-1}. \phi(p)^{-1}
  \exp(\int_{r' \in T(s(p))} \log(\psi(p,r')) d\mu)^{-1} \\
& = & \phi(p.q).  \phi(q)^{-1}.\phi(p)^{-1}. E(p,q) = \psi(p,q). E(p,q) \
,
\end{array}
\end{equation}
where

\begin{equation}
\begin{array}{lll}
E(p,q) & =  & Ad_{\psi(p,q)^{-1}}\exp(\int_{r \in T(s(q))} \log(\psi(p.q,r))
d\mu). \\
& &   Ad_{\phi(p)}\exp(\int_{r \in T(s(q))}^{-1} \log(\psi(q,r)) d\mu)^{-1}.
\exp(\int_{r' \in T(s(p))} \log(\psi(p,r')) d\mu)^{-1} \\
& = & \exp(\int_{r \in T(s(q))} \log(\psi(p,q)^{-1}.\psi(p.q,r).\psi(p,q)) d\mu). \\
& & \exp(\int_{r \in T(s(q))} \log(\phi(p).\psi(q,r)^{-1}.\phi(p)^{-1}) d\mu). \\
& & \exp(\int_{r \in T(s(q))} \log(\psi(p,q.r)^{-1}) d\mu) \ \
({\rm we \ changed} \ r' \ {\rm by} \ r = q^{-1}.r') \\
& = & \exp(\int_{r \in T} \log(A_1) d\mu). \exp(\int_{r \in T} \log(A_2) d\mu).
\exp(\int_{r \in T} \log(A_3) d\mu) \ ,
\end{array}
\end{equation}
where $T = T(s(q))$ and

\begin{equation}
\label{eqn:A1A2A3}
\begin{array}{l}
A_1 = \psi(p,q)^{-1}.\psi(p.q,r).\psi(p,q) \ , \\
A_2 = \phi(p).\psi(q,r)^{-1}.\phi(p)^{-1} \ , \\
A_3 = \psi(p,q.r)^{-1} \ .
\end{array}
\end{equation}

One verifies directly that

\begin{equation}
A_1.A_2.A_3 = \psi(p.q)^{-1} \ .
\end{equation}

Consider $A_1,A_2,A_3$ as maps from $M_3$ to $G$. By definition, $\Delta(\phi) =
\|\psi\|_0$. The inequality $\Delta(\phi) \leq 1$ in the hypothesis of Lemma
\ref{lem:C0-estimate}, together with the fact that the metric on $G$ is
bi-invariant, implies that

\begin{equation}
\label{eqn:AAA} \|A_1\|_0 = \|A_2\|_0 = \|A_3\|_0 = \|\psi\|_0 = \Delta(\phi) \leq 1
\end{equation}

Applying Inequalities (\ref{eqn:C0-exp}), (\ref{eqn:Ck_integral}) and
(\ref{eqn:AAA}) several times to $E(p,q)$, we get:

\begin{equation}
\begin{array}{lll}
\log E(p,q) & =  & \epsilon_1 + \int_{r \in T} \log(A_1) d\mu + \int_{r \in T}
\log(A_2) d\mu +
 \int_{r \in T} \log(A_3) d\mu \\
& = & \epsilon_1+  \int_{r \in T} [\log(A_1)+ \log(A_2) + \log(A_3)] d\mu \\
& = & \epsilon_1 +\epsilon_2 + \int_{r \in T} \log(A_1A_2A_3) d\mu \\
& = & \epsilon_1 + \epsilon_2 - \log (\psi(p,q))
\end{array}
\end{equation}
where $\epsilon_1$ and $\epsilon_2$ are some functions such that

\begin{equation}
\|\epsilon_1\|_0, \|\epsilon_2\|_0 \preceq \Delta(\phi)^2
\end{equation}

In other words, we have $\| \log (\psi(p,q)) + \log E(p,q)\|_0 \preceq
\Delta(\phi)^2$, which implies, by Inequality (\ref{eqn:C0-exp2}), that
$\|\psi.E\|_0 \preceq \Delta(\phi)^2$. But we have $\widehat{\psi}(p,q) = \psi(p,q).
E(p,q)$, therefore

\begin{equation}
\Delta(\widehat{\phi}) = \|\widehat{\psi}\|_0 =\|\psi.E\| \preceq \Delta(\phi)^2
\end{equation}
\hfill $\square$

Lemma \ref{lem:C0-estimate} immediately implies the uniform convergence (i.e.
convergence in $C^0$ topology) of the sequence of maps $\phi_n: M \to G$, defined
iteratively by $\phi_{n+1} = \widehat{\phi_n}$, beginning with an arbitrary smooth
map $\phi_1$ which satisfies the inequality $\Delta(\phi_1) \leq C_0/4$. (This
inequality can always be achieved by replacing $B$ by a sufficiently small invariant
neighborhood of $x_0$ in $B$ if necessary). Indeed, since $\|\psi_2\|_0 =
\Delta(\phi_2) \leq (\Delta(\phi_1))^2/C_0 \leq \Delta(\phi_1) \leq C_0/4 \leq 1/4$
by Lemma \ref{lem:C0-estimate}, where
\begin{equation}
\psi_n(p,q) = \phi_n(p.q).\phi_n(q)^{-1}. \phi_n(p)^{-1} \ ,
\end{equation}
we can define $\phi_2 = \widehat{\phi_1}$, and so on, hence $\phi_n$ is well defined
for all $n \in \bbN$. By recurrence on $n$, one can show easily  that we have
\begin{equation}
\label{eqn:C0_est} \|\psi_n\|_{0} \leq C_{0}. (b_0)^{2^{n}}  \ \ \forall n \in \bbN,
\ \ {\rm where} \ \ b_0 = \frac{1}{2} < 1 \ ,
\end{equation}
which implies in particular that $\sum_{n=1}^\infty\|\psi_n\|_0 < \infty$
(this is a very fast converging series). Put
\begin{equation}
\label{eqn:Psi_n}
\Psi_n(p,q) = \exp( \int_{q \in T(s(p))} \log(\psi_n(p.q)) d\mu) \ .
\end{equation}

Then $\|\Psi_n\|_0 \preceq \|\psi_n\|_0$ (by Inequalities (\ref{eqn:composition_est}) and
(\ref{eqn:Ck_integral})), which together with $\sum_{n=1}^\infty\|\psi_n\|_0 < \infty$
implies  that
\begin{equation}
\sum_{n=1}^\infty \|\Psi_n\|_0 < \infty \ .
\end{equation}

This last inequality implies the convergence of the product
$\Psi_n.\Psi_{n-1}...\Psi_1$ in $C^0$-topology when $n \to \infty$. But

\begin{equation}
\Psi_n.\Psi_{n-1}...\Psi_1 = \phi_{n+1}.\phi_1^{-1}
\end{equation}

Thus $\phi_n$ converges in $C^0$-topology when $n \to \infty$. Denote by
$\phi_\infty$ the limit
\begin{equation}
\phi_\infty = \lim_{n \to \infty} \phi_n
\end{equation}
Then $\phi_\infty$ is a continuous homomorphism from $M$ to $G$. It is also clear
that the restriction of $\phi_\infty$ to $G$ is the identity map from $G$ to itself.

It remains to show that $\phi_\infty$ is smooth. This is the purpose of the next
subsection, where we will show that for any $k \in \bbN$, $k \leq h$ if $M$ belongs
to the class $C^h$ only, we have $\phi_\infty = \lim_{n \to \infty} \phi_n$ in
$C^k$-topology as well.

\begin{remark} When $G$ is Abelian then the above computations show that
$\widehat{\psi} = 1_G$, i.e. $\widehat{\phi}$ is already a homomorphism
from $M$ to $G$, and there is no need to iterate our averaging process.
\end{remark}

\subsection{$C^k$ estimates}

Roughly speaking, we want to make estimations on $\psi_n$ in order to show that, if
$k$ does not exceed the smoothness class of the groupoid $M \rightrightarrows B$, then
$\sum_{n=1}^\infty \|\psi_n\|_k < \infty$. If this series converges, then similarly to
the previous subsection, we also have $\sum_{n=1}^\infty \|\Psi_n\|_k < \infty$ where
$\Psi_n = \phi_{n+1}.\phi_n^{-1}$ is given by formula (\ref{eqn:Psi_n}), hence the
product $\Psi_n.\Psi_{n-1}...\Psi_1$ converges in $C^k$-topology when $n \to \infty$,
implying that $\phi_n \to \phi_\infty$ in $C^k$-topology.

\begin{lem}
\label{lem:Ck-estimate} Let $k \in \bbN$ be a natural number which does not exceed
the smoothness class of the groupoid $M \rightrightarrows B$. Assume that
$\|\psi\|_0 = \Delta(\phi) \leq 1$ and $\|\phi\|_{k-1} \preceq 1$. Then we have:
\begin{equation}
\|\widehat{\psi}\|_k \preceq \|\psi\|_0\|\psi\|_k + \|\psi\|_{k-1} +
\|\psi\|_0\|\phi\|_{k-1} + \|\psi\|_0^2\|\phi\|_k
\end{equation}
\end{lem}

{\it Proof}. Assume that $\|\phi\|_{k-1} \preceq 1$ by hypothesis of Lemma
\ref{lem:Ck-estimate}. Then by Inequality (\ref{eqn:Ck_product}) and Inequality
(\ref{eqn:composition_est}), we have $\|\psi\|_{k-1} \preceq \|\phi\|_{k-1} \preceq
1$. Let $A_1,A_2,A_3$ be the functions defined by Equation (\ref{eqn:A1A2A3}). We
want to estimate them. For $A_3 = \psi(p,q.r^{-1})^{-1}$, using Inequality
(\ref{eqn:composition_est}), we get:

\begin{equation}
\|A_3\|_{k-1} \preceq \|\psi\|_{k-1} \preceq 1 \ \ {\rm and} \ \ \|A_3\|_{k} \preceq
\|\psi\|_{k} \ .
\end{equation}

For $A_1 =  \psi(p,q)^{-1}.\psi(p.q,r).\psi(p,q)$, using Inequality
(\ref{eqn:Ck_product}) (and Inequality (\ref{eqn:composition_est})), we also get

\begin{equation}
\|A_1\|_{k-1} \preceq \|\psi\|_{k-1} \preceq 1 \ \ {\rm and} \ \ \|A_1\|_{k} \preceq
\|\psi\|_{k} \ .
\end{equation}

The estimation of $A_2  = \phi(p).\psi(q,r)^{-1}.\phi(p)^{-1}$ is more complicated,
because it involves the function $\phi$ directly. Using Inequality
(\ref{eqn:Ck_conjugate}) we get
\begin{equation}
\|A_2\|_{k-1} \preceq \|\psi\|_{k-1} + \|\phi\|_{k-1}\|\psi\|_0 \preceq 1
\end{equation}
and
\begin{equation}
\|A_2\|_{k} \preceq \|\psi\|_{k} + \|\phi\|_{k}\|\psi\|_0 \ .
\end{equation}

Applying Inequality (\ref{eqn:Ck-exp}) and the above inequalities to $E(p,q)$, we
get that

\begin{equation}
\|\epsilon_1\|_k, \|\epsilon_2\|_k \preceq (\|\psi\|_k +\|\phi\|_k\|\psi\|_0)
\|\psi\|_0 + (\|\psi\|_{k-1} + \|\phi\|_{k-1}\|\psi\|_0)
\end{equation}

Moreover, we have

\begin{equation}
\|E\|_i \preceq \|\psi\|_i + \|\psi\|_0\|\phi\|_i \ \ \forall i=0,...,k.
\end{equation}

Now applying Inequality (\ref{eqn:Ck-exp2}) and the last two inequalities, we get

\begin{equation}
\begin{array}{lll}
\|\widehat{\psi}\|_k & = & \|\exp(\log(\psi)).\exp(\log(E)) \|_k \\
& \preceq & \| \log(\psi) + \log(E)\|_k + \|\psi\|_0\|E\|_k + \|\psi\|_k\|E\|_0 +
 \|\psi\|_{k-1} + \|E\|_{k-1} \\
& = & \|\epsilon_1 + \epsilon_2\|_k + \|\psi\|_0\|E\|_k + \|\psi\|_k\|E\|_0 +
 \|\psi\|_{k-1} + \|E\|_{k-1} \\
& \preceq & \|\psi\|_0\|\psi\|_k + \|\psi\|_{k-1} + \|\psi\|_0\|\phi\|_{k-1} +
\|\psi\|_0^2\|\phi\|_k
\end{array}
\end{equation}
\hfill $\square$

\begin{lem}
\label{lem:Ck-estimate2}
With the assumptions of Lemma \ref{lem:Ck-estimate} we have:
\begin{equation}
\|\widehat{\phi}\|_k - \|\phi\|_k \preceq \|\psi\|_k + \|\psi\|_0 \|\phi\|_k
\end{equation}
\end{lem}
{\it Proof}. Applying Inequality (\ref{eqn:Ck_product2}) to $\widehat{\phi} = \Psi.\phi$,
we get
\begin{equation}
\|\widehat{\phi}\|_k - \|\phi\|_k \preceq \|\Psi\|_k + \|\Psi\|_0 \|\phi\|_k
\end{equation}
Now replace $\|\Psi\|_0$ by $\|\psi\|_0$ and $\|\Psi\|_k$ by $\|\psi\|_k$.
\hfill $\square$

\begin{lem}
\label{lem:Ck-estimate3}
Assume that $\phi_1$ is a map from $M$ to $G$ such that
$\Delta(\phi_1) < C_0/4$, and that $\phi_{n+1} = \widehat{\phi_n}$ for any $n \in \bbN$,
as in the previous subsection. Let $k$ be
a natural number which does not exceed the smoothness class of the groupoid
$M \rightrightarrows B$. Then there is a finite positive number $D_k > 0$ and
a positive number $0 < b_k < 1$, such that for any $n \in \bbN$ the
following two inequalities hold:
\begin{equation}
\label{eqn:Ck_phi}
\|\phi_n\|_k \leq D_k.(1 - 2^{-n})
\end{equation}
and
\begin{equation}
\label{eqn:Ck_psi}
\|\psi_n\|_k \leq D_k.(b_k)^{2^n} \ .
\end{equation}
\end{lem}

\emph{Proof of Theorem \ref{thm:LPG}}. Remark that if Inequality
(\ref{eqn:Ck_phi}) is satisfied, then in particular the $k$-norms
$\|\phi_n\|_k$ ($n \in \bbN$) is bounded by $D_k$, and this is enough to
imply (by Ascoli theorem) that $\phi_\infty = \lim_{n \to \infty} \phi_n$
is of class $C^k$. Inequality (\ref{eqn:Ck_psi}) is also a sufficient
condition for the $C^k$-smoothness of $\phi_\infty$, because it implies in
particular that $\sum_{n=1}^\infty \|\Psi_n\|_k \preceq \sum_{n=1}^\infty
\|\psi_n\|_k < \infty$, which in turns implies that the sequence of maps
$(\phi_n)$ converges in $C^k$-topology. Thus the smoothness of the
homomorphism $\phi_\infty: M \to G$, and hence Theorem \ref{thm:LPG},
follows directly from the above Lemma, and in fact we need only one of the
above two inequalities (\ref{eqn:Ck_psi}) and (\ref{eqn:Ck_phi}) in order
to prove our main result. But these two inequalities are closely related
(they need each other in a proof by induction), so we put them together.
\hfill $\square$

\emph{Proof of Lemma \ref{lem:Ck-estimate3}}. We will prove it by induction on $k$.
When $k=0$, Lemma \ref{lem:Ck-estimate3} is already proved in the previous section
(with $b_0 = 1/2$). Let us now assume that Inequalities (\ref{eqn:Ck_psi}) and
(\ref{eqn:Ck_phi}) are true at the level $k-1$ (i.e. if we replace $k$ by $k-1$). We
will show that they are true at the level $k$.

We will choose an (arbitrary) number $b_k > 0$ such that $1 > b_k > b_k^2 > b_{k-1}, b_0$.
(For example, one can put $b_0 = 1/2$ and then $b_k = (b_{k-1})^{1/3}$ by recurrence).
What will be important for us is that $b_0/b_k,b_{k-1}/b_k^2$ and $b_0/b_{k}^2$ are
positive numbers which are strictly smaller than 1.

It follows from Lemma \ref{lem:Ck-estimate} and Lemma \ref{lem:Ck-estimate2} that
there exist two positive numbers $c_1$ and $c_2$ (which do not depend on $n$) such
that we have, for any $n \in \bbN$:
\begin{equation}
\label{eqn:psi_est2} \|\psi_{n+1}\|_k \leq c_1(\|\psi_n\|_0\|\psi_n\|_k +
\|\psi_n\|_{k-1} + \|\psi_n\|_0\|\phi_n\|_{k-1} + \|\psi_n\|_0^2\|\phi_n\|_k)
\end{equation}
and
\begin{equation}
\|\phi_{n+1}\|_k - \|\phi_n\|_k \leq c_2(\|\psi_n\|_k + \|\psi_n\|_0 \|\phi_n\|_k) \
.
\end{equation}

We will now prove Inequalities (\ref{eqn:Ck_psi}) and
(\ref{eqn:Ck_phi}) by induction on $n$. There exists a natural number $n_0$ such that
for any $n > 0$ we have

\begin{equation}
Q_1 := D_0 \left(\frac{b_0}{b_k}\right)^{2^n} +  \left(\frac{b_{k-1}}{b_k^2}\right)^{2^n} +
D_0 \left(\frac{b_0}{b_k^2}\right)^{2^n}
+ D_0^2 \left(\frac{b_0}{b_k}\right)^{2^{n+1}} \leq \frac{1}{c_1}
\end{equation}
and
\begin{equation}
Q_2 := (b_k)^{2^n} + D_0 (b_0)^{2^n} \leq \frac{2^{-n-1}}{ c_2} \ .
\end{equation}

By choosing $D_k$ large enough, we can assume that Inequalities (\ref{eqn:Ck_psi}) and
(\ref{eqn:Ck_phi}) are satisfied for any $n \leq n_0$. We will also assume that
$D_k \geq D_{k-1}$. let us now show that if Inequalities (\ref{eqn:Ck_psi}) and
(\ref{eqn:Ck_phi}) are satisfied for some $n \geq n_0$ then they are still satisfied
when we replace $n$ by $n+1$. (This is the last step in our induction process).

Indeed, for $\|\psi_{n+1}\|_k$, using Inequality (\ref{eqn:psi_est2}) and the induction
hypothesis, we get

\begin{equation}
\begin{array}{ll}
\|\psi_{n+1}\|_k & \leq c_1\left(\|\psi_n\|_k \|\psi_n\|_0 + \|\psi_n\|_{k-1} +
\|\psi_n\|_0 \|\phi_n\|_{k-1} + \|\psi_n\|_0^2 \|\phi_n\|_k \right) \\
& \leq c_1 \left( D_k(b_k)^{2^n}D_0(b_0)^{2^n} + D_{k-1}(b_{k-1})^{2^n} + D_0 (b_0)^{2^n}
D_{k-1} + D_0^2 (b_0)^{2^{n+1}}D_k \right) \\
& \leq D_k c_1 \left( (b_k)^{2^n}D_0(b_0)^{2^n} + (b_{k-1})^{2^n} + D_0 (b_0)^{2^n}
 + D_0^2 (b_0)^{2^{n+1}} \right) \\
 & = D_k c_1 Q_1 (b_k)^{2^{n+1}} \\
 & \leq D_k (b_k)^{2^{n+1}}
\end{array}
\end{equation}

Similarly, for $\|\phi_{n+1}\|_k$ we have:

\begin{equation}
\begin{array}{ll}
\|\phi_{n+1}\|_k & \leq \|\phi_n\|_k + c_2\left(\|\psi_n\|_k + \|\psi_n\|_0\|\phi_n\|_{k}
\right) \\
& \leq D_{k}(1 - 2^{-n}) + c_2 \left( D_k(b_k)^{2^n} + D_0(b_0)^{2^n}.D_k \right) \\
& \leq D_{k}(1 - 2^{-n}) + D_k c_2 [(b_k)^{2^n} + D_0(b_0)^{2^n}] \\
 & = D_{k}(1 - 2^{-n})  + D_kc_2 Q_2 \\
 & \leq D_k (1 - 2^{-n}) + D_k 2^{-n-1} = D_k(1 - 2^{-n-1})
\end{array}
\end{equation}

\hfill $\square$

\begin{remark} If we start with a near-homomorphism from $M$ to a
compact Lie group $H$ different from $G$, then our iterative averaging
method still yields a homomorphism from $M$ to $H$. Thus, in a sense, our
result generalizes a result of Grove, Karcher and Ruh
\cite{GrKaRu-Group1974} on near-homomorphisms between compact groups.
\end{remark}




\vspace{0.3cm} {\bf Acknowledgements}. I would like to thank Alan
Weinstein for getting my attention to the problem, and for many
stimulating discussions. Thanks to Philippe Monnier for finding some
errors in the first version of this paper.

\bibliographystyle{amsplain}
\providecommand{\bysame}{\leavevmode\hbox to3em{\hrulefill}\thinspace}

\end{document}